\documentclass[12pt]{article} 
\usepackage{amsmath}
\usepackage{amssymb} 
\usepackage{amscd}

\usepackage{theorem} 

\usepackage{xspace}

\def\im{{\mbox{Im}}}

\def\ker{{\mbox{Ker}}}
\def\Hom {{\mbox{Hom}}}

\def\End{{\mbox{End}}}

\def\cala{{\cal A}} 
\def\calb{{\cal B}}

\def\calc{{\cal C}}

\def\bbbone{\mbox{\rm 1\hspace {-.6em} l}}

\def\hg{{\mathbf H}}
\def\alg{{\mathbf{Alg}}} 
\def\vect{{\mathbf{Vect}}} 
\def\Homg {{\mathbf{Hom}}}
\def\homg {{\mathbf{hom}}}
\def\endg{{\mathbf{end}}}

\newtheorem{theorem}{THEOREM}
\newtheorem{lemma}{LEMMA} 
\newtheorem{proposition}{PROPOSITION}
\newtheorem{corollary}{COROLLARY}

\begin{document}

\baselineskip=0.7cm

\begin{center} 
 \thispagestyle{empty}
{\large\bf HOMOGENEOUS ALGEBRAS}
\end{center} 
\vspace{0.75cm}

\begin{center} 
 Roland BERGER
 \footnote{LARAL, Facult\'e des Sciences et Techniques,
 23 rue P. Michelon,
 F-42023 Saint-Etienne Cedex 2, France\\
Roland.Berger@univ-st-etienne.fr\\
 },
 Michel DUBOIS-VIOLETTE
\footnote{Laboratoire de Physique Th\'eorique, UMR 8627, Universit\'e Paris XI,
B\^atiment 210, F-91 405 Orsay Cedex, France\\
Michel.Dubois-Violette$@$th.u-psud.fr\\
},
Marc WAMBST
\footnote{Institut de Recherche Math\'ematique Avanc\'ee,
Universit\'e Louis Pasteur - C.N.R.S.,\\
7 rue Ren\'e Descartes,
F-67084 Strasbourg Cedex, France\\
wambst@math.u-strasbg.fr\\
}
\end{center} \vspace{1cm}

\begin{center} \today \end{center}

\vspace {1cm}

\begin{abstract}
 Various concepts associated with quadratic algebras admit natural 
 generalizations when the quadratic algebras are replaced by graded 
 algebras which are finitely generated in degree 1 with homogeneous 
 relations of degree $N$. Such algebras are referred to as {\sl 
 homogeneous algebras of degree $N$}. In particular it is shown that 
 the Koszul complexes of quadratic algebras generalize as 
 $N$-complexes for homogeneous algebras of degree $N$.

\end{abstract}

  \vspace{1,5cm}
\noindent LPT-ORSAY 02-08

\newpage

\section{Introduction and Preliminaries}

Our aim is to generalize the various concepts associated with 
quadratic algebras as described in \cite{YuM2} when the quadratic 
algebras are replaced by the homogeneous algebras of degree $N$ with 
$N\geq 2$  ($N=2$ is the case of quadratic algebras). Since the 
generalization is natural and relatively straightforward, the 
treatment of \cite{YuM1}, \cite{YuM2} and \cite{JLL} will be directly 
adapted to homogeneous algebras of degree $N$. In other words we 
dispense ourselves to give a review of the case of quadratic algebras 
(i.e. the case $N=2$) by referring to the above quoted nice 
treatments. In proceeding to this adaptation, we shall make use of the 
following slight elaboration of an ingredient of the elegant presentation of \cite{JLL}.

\begin{lemma}
 Let $A$ be an associative algebra with product denoted by $m$, let 
 $C$ be a coassociative coalgebra with coproduct denoted by $\Delta$ 
 and let $\Hom_{\mathbb K}(C,A)$ be equipped with its structure of 
 associative algebra for the convolution product 
 $(\alpha,\beta)\mapsto \alpha\ast \beta=m\circ (\alpha\otimes 
 \beta)\circ \Delta$. Then one defines an algebra-homomorphism 
 $\alpha\mapsto d_{\alpha}$ of $\Hom_{\mathbb K}(C,A)$ into the 
 algebra $\End_{A}(A\otimes C)=\Hom_{A}(A\otimes C, 
 A\otimes C)$ of endomorphisms of the left $A$-module $A\otimes C$ by defining 
 $d_{\alpha}$ as the composite 
 \[
 A\otimes C \stackrel{I_{A} \otimes \Delta}{\longrightarrow} A\otimes 
 C\otimes C\stackrel{ I_{A}\otimes \alpha \otimes 
I_{C}}{\longrightarrow}  A\otimes A\otimes C \stackrel{m\otimes I_{C}
}{\longrightarrow} A\otimes C
 \]
 for $\alpha\in \Hom_{\mathbb K}(C,A)$.
 \end{lemma}
 The proof is straightforward, $d_{\alpha}\circ 
 d_{\beta}=d_{\alpha\ast\beta}$ follows easily from the 
 coassociativity of $\Delta$ and the associativity of $m$. As pointed 
 out in \cite{JLL} one obtains a graphical version (``electronic 
 version'') of the proof by using the usual graphical version of the 
 coassociativity of $\Delta$ combined with the usual graphical 
 version of the associativity of $m$. The left $A$-linearity of $d_\alpha$ is straightforward.\\
 
 In the above statement as well as in the following, all vector 
 spaces, algebras, coalgebras are over a fixed field $\mathbb K$. Furthermore unless otherwise specified the 
 algebras are unital associative and the coalgebras are counital 
 coassociative. For instance in the previous case, if $\bbbone$ is the 
 unit of $A$ and $\varepsilon$ is the counit of $C$, then the unit 
 of $\Hom_{\mathbb K}(C,A)$ is the linear mapping $\alpha\mapsto 
 \varepsilon(\alpha)\bbbone$ of $C$ into $A$. In Lemma 1 the left $A$-module structure on $A\otimes C$ is the obvious one given by 
\[
x(a\otimes c) = (xa)\otimes c
\]
for any $x\in A$, $a\in A$ and $c\in C$.\\

Besides the fact that it is natural to generalize for other degrees what exists for quadratic algebras, this paper produces a very natural class of $N$-complexes which generalize the Koszul complexes of quadratic algebras \cite{YuM1}, \cite{YuM2}, \cite{Wam1}, \cite{JLL}, \cite{AG} and which are not of simplicial type. By $N$-complexes of simplicial type we here mean $N$-complexes associated with simplicial modules and $N$-th roots of unity in a very general sense \cite{D-V2} which cover cases considered e.g. in \cite{May}, \cite{Kap}, \cite{D-VK}, \cite{D-V}, \cite{KW} the generalized homology of which has been shown to be equivalent to the ordinary homology of the corresponding simplicial modules \cite{D-V2}. This latter type of constructions and results has been recently generalized to the case of cyclic modules \cite{Wam3}. In spite of the fact that they compute the ordinary homology of the simplicial modules, the usefulness of these $N$-complexes of simplicial type comes from the fact that they can be combined with other $N$-complexes \cite{D-VT1}, \cite{D-VT2}. In fact the BRS-like construction \cite{BRS} of \cite{D-VT2} shows that spectral sequences arguments (e.g. in the form of a generalization of the homological perturbation theory \cite{Sta}) are still working for $N$-complexes. Other nontrivial classes of $N$-complexes which are not of simplicial type are the universal construction of \cite{D-VK} and the $N$-complexes of \cite{D-VH}, \cite{D-VH2} (see also in \cite{D-V4} for a review). It is worth noticing here that elements of homological algebra for $N$-complexes have been developed in \cite{KW} and that several results for $N$-complexes and more  generally $N$-differential modules like Lemma 1 of \cite{D-V2} have no nontrivial counterpart for ordinary complexes and differential modules. It is also worth noticing that besides the above mentioned examples,  various problems connected with theoretical physics implicitly involve exotic $N$-complexes (see e.g. \cite{RK1}, \cite{RK2}).\\

In the course of the paper we shall point out the possibility of generalizing the approach based on quadratic algebras of \cite{YuM2} to quantum spaces and quantum groups by replacing the quadratic algebras by $N$-homogeneous ones. Indeed one also has in this framework internal $\endg$, etc. with similar properties.\\

Finally we shall revisit in the present context the approach of \cite{RB3}, \cite{RB4} to Koszulity for $N$-homogeneous algebras. This is in order since as explained below, the generalization of the Koszul complexes introduced in this paper for $N$-homogeneous algebras is a canonical one. We shall explain why a definition based on the acyclicity of the $N$-complex generalizing the Koszul complex is inappropriate and we shall identify the ordinary complex introduced in \cite{RB3} (the acyclicity of which is the definition of Koszulity of \cite{RB3}) with a complex obtained by contraction from the above Koszul $N$-complex. Furthermore we shall show the uniqueness of this contracted complex among all other ones. Namely we shall show that the acyclicity of any other complex (distinct from the one of \cite{RB3}) obtained by contraction of the Koszul $N$-complex leads for $N\geq 3$ to an uninteresting (trivial) class of algebras.\\

Some examples of Koszul homogeneous algebras of degree $>2$ are given
in \cite{RB3}, including a certain cubic Artin-Schelter regular 
algebra \cite{AS}. Recall that Koszul quadratic
algebras arise in several topics as algebraic geometry \cite{Kempf},
representation theory \cite{BGS}, quantum groups \cite{YuM1},
\cite{YuM2}, \cite{Wam1}, \cite{Wam2}, Sklyanin algebras \cite{SmSt},
\cite{Tvdb}. A classification of the Koszul quadratic algebras with
two generators over the complex numbers is performed in \cite{RB2}. 
Koszulity of non-quadratic algebras and each of the above items
deserve further attention.\\

The plan of the paper is the following.\\
 In Section 2 we define the duality and the two (tensor) products which are exchanged by the duality for homogeneous algebras of degree $N$ ($N$-homogeneous algebras). These are the direct extension to arbitrary $N$ of the concepts defined for quadratic algebras ($N=2$), \cite{YuM1}, \cite{YuM2}, \cite{JLL} and our presentation here as well as in Section 3 follows closely the one of reference \cite{YuM2} for quadratic algebras.\\
In Section 3 we elaborate the categorical setting and we point out the conceptual reason for the occurrence of $N$-complexes in the framework of $N$-homogeneous algebras. We also sketch in this section a possible extension of the approach of \cite{YuM2} to quantum spaces and quantum groups in which relations of degree $N$ replace the quadratic ones.\\
In Section 4 we define the $N$-complexes which are the generalizations for homogeneous algebras of degree $N$ of the Koszul complexes of quadratic algebras \cite{YuM1}, \cite{YuM2}. The definition of the cochain $N$-complex $L(f)$ associated with a morphism $f$ of $N$-homogeneous algebras follows immediately from the structure of the unit object $\wedge_N\{d\}$ of one of the (tensor) products of $N$-homogeneous algebras. We give three equivalent definitions of the chain $N$-complexe $K(f)$: A first one by dualization of the definition of $L(f)$, a second one which is an adaptation of \cite{JLL} by using Lemma 1, and a third one which is a component-wise approach. It is pointed out in this section that one cannot generalize naively the notion of Koszulity for $N$-homogeneous algebras with $N\geq 3$ by the acyclicity of the appropriate Koszul $N$-complexes.\\
In Section 5, we recall the definition of Koszul homogeneous algebras of \cite{RB3} as well as some results of \cite{RB3}, \cite{RB4} which justify this definition. It is then shown that this definition of Koszulity for homogeneous $N$-algebras is optimal within the framework of the appropriate Koszul $N$-complex.\\

Let us give some indications on our notations. Throughout the paper the symbol $\otimes$ denotes the tensor product over the basic field $\mathbb K$. Concerning the generalized homology of $N$-complexes we shall use the notation of \cite{Kap} which is better adapted than other ones to the case of chain $N$-complexes, that is if $E=\oplus_n E_n$ is a chain $N$-complex with $N$-differential $d$, its generalized homology is denoted by $_pH(E)=\oplus_{n\in \mathbb Z}\> \> _pH_n(E)$
 with 
\[
_pH_n(E)=\ker (d^p:E_n\rightarrow E_{n-p})/\im (d^{N-p}:E_{n+N-p}\rightarrow E_n)
\]
for $p\in\{ 1,\dots, N-1\}$, ($n\in \mathbb Z$).

\section{Homogeneous algebras of degree N}

Let $N$ be an integer with $N\geq 2$. A {\sl homogeneous algebra of 
degree $N$} or {\sl $N$-homogeneous algebra} is an algebra of the form
\begin{equation}
\cala = A(E,R)=T(E)/(R)
\label{eq1}
\end{equation}
where $E$ is a finite-dimensional vector space (over $\mathbb K$), 
$T(E)$ is the tensor algebra of $E$ and $(R)$ is the two-sided ideal 
of $T(E)$ generated by a linear subspace $R$ of $E^{\otimes^N}$. The 
homogeneity of $(R)$ implies that $\cala$ is a graded algebra 
$\cala=\oplus_{n\in \mathbb N}\cala_n$ with 
$\cala_n=E^{\otimes^n}$ for $n<N$ and 
$\cala_n=E^{\otimes^n}/\sum_{r+s=n-N}E^{\otimes^r}\otimes R\otimes 
E^{\otimes^s}$ for $n\geq N$ where we have set $E^{\otimes^0}=\mathbb 
K$ as usual. Thus $\cala$ is a graded algebra which is connected 
($\cala_0=\mathbb K)$, generated in degree 1 $(\cala_1=E)$ with the 
ideal of relations among the elements of $\cala_1=E$ generated by 
$R\subset E^{\otimes^N}=(\cala_1)^{\otimes^N}$.\\

A {\sl morphism of $N$-homogeneous algebras} $f:A(E,R)\rightarrow 
A(E',R')$ is a linear mapping $f:E\rightarrow E'$ such that 
$f^{\otimes^N}(R)\subset R'$. Such a morphism is a homomorphism of 
unital graded algebras. Thus one has a category $\hg_{N}\alg$ of 
$N$-homogeneous algebras and the forgetful functor 
$\hg_{N}\alg\rightarrow \vect$, $\cala\mapsto E$, from $\hg_{N}\alg$ to the category 
$\vect$ of finite-dimensional vector spaces (over $\mathbb K$).\\

Let $\cala=A(E,R)$ be a $N$-homogeneous algebra. One defines {\sl its 
dual} $\cala^!$ to be the $N$-homogeneous algebra 
$\cala^!=A(E^\ast,R^\perp)$ where $E^\ast$ is the dual vector space 
of $E$ and where $R^\perp\subset 
E^{\ast\otimes^N}=(E^{\otimes^N})^\ast$ is the annihilator of $R$ 
i.e. the subspace $\{\omega \in (E^{\otimes^N})^\ast\vert 
\omega(x)=0,\ \  \forall x\in R\}$ of $(E^{\otimes^N})^\ast$ 
identified with $E^{\ast\otimes^N}$. One has canonically 
\begin{equation}
(\cala^!)^!=\cala
\label{eq2}
\end{equation}
and if $f:\cala\rightarrow \cala'=A(E',R')$, is a morphism of 
$\hg_{N}\alg$, the transposed of $f:E\rightarrow E'$ is a linear 
mapping of $E^{\prime\ast}$ into $E^\ast$ which induces the morphism 
$f^!:(\cala')^!\rightarrow \cala^!$ of $\hg_{N}\alg$ so 
$(\cala\mapsto \cala^!$, $f\mapsto f^!$) is a contravariant 
(involutive) functor.\\

Let $\cala=A(E,R)$ and $\cala'=A(E',R')$ be $N$-homogeneous 
algebras; one defines $\cala\circ \cala'$ and $\cala\bullet\cala'$ by 
setting
\[
 \cala\circ \cala'=A(E\otimes E',\ \  \pi_{N}(R\otimes 
 E^{\prime\otimes^N}+E^{\otimes^N}\otimes R'))
 \]
 \[\cala\bullet \cala'=A(E\otimes E',\ \  \pi_{N}(R\otimes R'))
 \]
 where $\pi_{N}$ is the permutation
 \begin{equation}
 (1,2,\dots,2N)\mapsto (1,N+1,2,N+2,\dots,k,N+k,\dots,N,2N)
 \label{eq3}
 \end{equation}
belonging to the symmetric group $S_{2N}$ acting as usually on the 
factors of the tensor products. One has canonically
\begin{equation}
(\cala\circ\cala')^!=\cala^!\bullet \cala^{\prime !},\ \ 
(\cala\bullet \cala')^!=\cala^!\circ \cala^{\prime !}
\label{eq4}
\end{equation}
which follows from the identity $\{ R\otimes E^{\prime 
\otimes^N}+E^{\otimes^N}\otimes R'\}^\perp = R^\perp \otimes 
R^{\prime\perp}$. On the other hand the inclusion $R\otimes R'\subset 
R\otimes E^{\prime\otimes^N}+E^{\otimes^N}\otimes R'$ induces an 
surjective algebra-homomorphism $p:\cala\bullet \cala'\rightarrow 
\cala\circ \cala'$ which is of course a morphism of $\hg_{N}\alg$.\\

It is worth noticing here that in contrast with what happens for 
quadratic algebras if $\cala$ and $\cala'$ are homogeneous algebras of 
degree $N$ with $N\geq 3$ then the tensor product algebra 
$\cala\otimes \cala'$ is no more a $N$-homogeneous algebra. 
Nevertheless there still exists an injective homomorphism of unital 
algebra $i:\cala\circ \cala' \rightarrow \cala\otimes \cala'$ doubling 
the degree which we now describe. Let $\tilde\imath:T(E\otimes E')\rightarrow T(E)\otimes T(E')$ be the 
injective linear mapping which restricts as
\[
\tilde\imath 
=\pi^{-1}_{n}:(E\otimes E')^{\otimes^n}\rightarrow 
E^{\otimes^n}\otimes E^{\prime\otimes^n}
\]
on $T^n(E\otimes E')=(E\otimes E')^{\otimes^n}$ for any $n\in\mathbb 
N$. It is straightforward that $\tilde\imath$ is an 
algebra-homomorphism which is an isomorphism onto the subalgebra 
$\oplus_{n}E^{\otimes^n}\otimes E^{\prime\otimes^n}$ of $T(E)\otimes 
T(E')$. The following proposition is not hard to verify.

\begin{proposition}
 Let $\cala=A(E,R)$ and $\cala'=A(E',R')$ be two $N$-homoge\-neous 
 algebras. Then $\tilde\imath$ passes to the quotient and induces an 
 injective homomorphism $i$ of unital algebras of $\cala\circ\cala'$ 
 into $\cala\otimes \cala'$. The image of $i$ is the subalgebra 
 $\oplus_{n}\cala_n\otimes\cala^{\prime}_n$ of 
 $\cala\otimes\cala^\prime$.
 \end{proposition}
 The proof is almost the same as for quadratic algebras \cite{YuM2}.\\
 
 \noindent \underbar{Remark}. As pointed out in \cite{YuM2}, any 
 finitely related and finitely generated graded algebra (so in 
 particular any $N$-homogeneous algebra) gives rise to a quadratic 
 algebra. Indeed if $\cala=\oplus_{n\geq 0}\cala_n$ is a graded 
 algebra, define $\cala^{(d)}$ by setting 
 $\cala^{(d)}=\oplus_{n\geq 0}\cala_{nd}$. Then it was 
 shown in \cite{BF} that if $\cala$ is generated by the 
 finite-dimensional subspace $\cala_1$ of its elements of degree 1 
 with the ideal of relations generated by its components of degree 
 $\leq r$, then the same is true for $\cala^{(d)}$ with $r$ replaced 
 by $2+(r-2)/d$.
 
 \section{Categorical properties}
 
 Our aim in this section is to investigate the properties of the 
 category $\hg_{N}\alg$. We follow again closely \cite{YuM2} replacing 
 the quadratic algebras considered there by the $N$-homogeneous 
 algebras.\\
 
  Let $\cala=A(E,R)$, $\cala'=A(E',R')$ and $\cala''=A(E'',R'')$ be 
 three homogeneous algebras of degree $N$. Then the isomorphisms 
 $E\otimes E'\simeq E'\otimes E$ and $(E\otimes E')\otimes E''\simeq 
 E\otimes (E'\otimes E'')$ of $\vect$ induce corresponding 
 isomorphisms $\cala\circ \cala'\simeq \cala'\circ\cala$ and 
 $(\cala\circ\cala')\circ\cala''\simeq 
 \cala\circ(\cala'\circ\cala'')$ of $N$~-~homogeneous algebras (i.e. of $\hg_{N}\alg$). Thus $\hg_{N}\alg$ endowed with $\circ$ is a 
 tensor category \cite{DM} and furthermore to the 1-dimensional 
 vector space $\mathbb K t\in \vect$ which is a unit object of 
 $(\vect,\otimes)$ corresponds the polynomial algebra $\mathbb 
 K[t]=A(\mathbb K t,0)\simeq T(\mathbb K)$ as unit object of 
 ($\hg_{N}\alg,\circ)$. In fact the isomorphisms $\mathbb K[t]\circ 
 \cala\simeq \cala\simeq \cala\circ \mathbb K[t]$ are obvious in 
 $\hg_{N}\alg$. Thus one has Part $(i)$ of the following theorem.
 
 \begin{theorem}
 The category $\hg_{N}\alg$ of $N$-homogeneous algebras has the 
 following properties $(i)$ and $(ii)$\\
 $(i)$ $\hg_{N}\alg$ endowed with $\circ$ is a tensor category with unit 
 object $\mathbb K[t]$.\\
 $(ii)$ $\hg_{N}\alg$ endowed with $\bullet$ is a tensor category with 
 unit object $\wedge_{N}\{d\}=\mathbb K[t]^!$.
 \end{theorem}
 Part $(ii)$ follows from $(i)$ by the duality $\cala\mapsto \cala^!$. 
 In fact $(i)$ and $(ii)$ are equivalent in view of (2) and (4).\\
 
 The $N$-homogeneous algebra $\wedge_{N}\{d\}=\mathbb K[t]^!\simeq 
 T(\mathbb K)/\mathbb K^{\otimes^N}$ is the (unital) graded algebra 
 generated in degree one by $d$ with relation $d^N=0$. Part $(ii)$  of 
 Theorem 1 is the very reason for the appearance of $N$-complexes in 
 the present context, remembering the obvious fact that graded 
 $\wedge_{N}\{d\}$-module and $N$-complexe are the same thing.
 
 \begin{theorem}
  The functorial isomorphism in $\vect$
  \[
  \Hom_{\mathbb K}(E\otimes E',E'')=\Hom_{\mathbb K}(E,E^{\prime 
  \ast}\otimes E'')
  \]
  induces a corresponding functorial isomorphism
  \[
  \Hom (\cala\bullet\calb, \calc)=\Hom(\cala,\calb^!\circ \calc)
  \]
  in $\hg_{N}\alg$, $($setting $\cala=A(E,R)$, $\calb=A(E',R')$ and 
  $\calc=A(E'',R''))$.
  \end{theorem}
  Again the proof is the same as for quadratic algebras \cite{YuM2}. 
  It follows that the tensor category $(\hg_{N}\alg,\bullet)$ has an 
  internal $\Homg$ \cite{DM} given by
  \begin{equation}
  \Homg (\calb, \calc)= \calb^!\circ \calc
  \label{eq5}
  \end{equation}
  for two $N$-homogeneous algebras $\calb$ and $\calc$. Setting 
  $\cala=A(E,R)$, $\calb=A(E',R')$ and $\calc=A(E'',R'')$ one verifies 
  that the canonical linear mappings $(E^\ast\otimes E')\otimes 
  E\rightarrow E'$ and $(E^{\prime\ast}\otimes E'')\otimes (E^\ast\otimes 
  E')\rightarrow E^\ast\otimes E''$ induce products
  \begin{equation}
  \mu: \Homg(\cala,\calb)\bullet \cala \rightarrow \calb
  \label{eq6}
  \end{equation}
  \begin{equation}
  m : \Homg(\calb, \calc)\bullet\Homg (\cala,\calb)\rightarrow 
  \Homg(\cala,\calc)
  \label{eq7}
  \end{equation}
  these internal products as well as their associativity properties 
  follow more generally from the formalism of tensor categories 
  \cite{DM}.\\
  
  Following \cite{YuM2}, define 
 $\homg (\cala,\calb)=\Homg(\cala^!,\calb^!)^!=\cala^!\bullet\calb$. 
  Then one obtains by duality from (\ref{eq6}) and (\ref{eq7}) 
  morphisms
  \begin{equation}
  \delta_{\circ}  : \calb \rightarrow \homg (\cala,\calb)\circ \cala
  \label{eq8}
  \end{equation}
  \begin{equation}
  \Delta_{\circ}  :\homg (\cala,\calc)\rightarrow \homg (\calb, \calc)\circ 
  \homg (\cala,\calb)
  \label{eq9}
  \end{equation}
  satisfying the corresponding coassociativity properties from which 
  one obtains by composition with the corresponding homomorphisms $i$ 
  the algebra homomorphisms
  
  \begin{equation}
  \delta  : \calb \rightarrow \homg (\cala,\calb)\otimes \cala
  \label{eq10}
  \end{equation}
  \begin{equation}
  \Delta  :\homg (\cala,\calc)\rightarrow \homg (\calb, \calc)\otimes
  \homg (\cala,\calb)
  \label{eq11}
  \end{equation}

  \begin{theorem}
   Let $\cala=A(E,R)$ be a $N$-homogeneous algebra. Then the 
   ($N$-homogeneous) algebra $\endg(\cala)=\cala^!\bullet \cala=\homg 
   (\cala,\cala)$ endowed with the coproduct $\Delta$ becomes a 
   bialgebra with counit $\varepsilon: \cala^!\bullet 
   \cala\rightarrow \mathbb K$ induced by the duality 
   $\varepsilon=\langle \cdot,\cdot\rangle:E^\ast\otimes 
   E\rightarrow \mathbb K$ and $\delta$ defines on $\cala$ a structure of left 
   $\endg(\cala)$-comodule.
   \end{theorem}
   
   \section{The $N$-complexes $L(f)$ and  $K(f)$}
   
   Let us apply Theorem 2 with $\cala=\wedge_{N}\{d\}$ and use 
   Theorem 1 $(ii)$. One has
   \begin{equation}
 \Hom(\calb,\calc)=\Hom(\wedge_{N}\{d\}, \calb^!\circ \calc)
 \label{eq12}
 \end{equation}
 and we denote by $\xi_{f}\in \calb^!\circ\calc$ the image of $d$ 
 corresponding to the morphism $f\in \Hom(\calb,\calc)$. One has 
 $(\xi_{f})^N=0$ and by using the injective algebra-homomorphism 
 $i:\calb^!\circ \calc\rightarrow \calb^!\otimes \calc$ of Proposition 1 
 we let $d$ be the left multiplication by $i(\xi_{f})$ in 
 $\calb^!\otimes \calc$. One has $d^N=0$ so, equipped with the 
 appropriate graduation, $(\calb^!\otimes \calc, d)$ 
 is a $N$-complex which will be denoted by $L(f)$. In the case where 
 $\cala=\calb=\calc$ and where $f$ is the identity mapping $I_{\cala}$ 
 of $\cala$ onto itself, this $N$-complex will be denoted by $L(\cala)$. 
 These $N$-complexes are the generalizations of the Koszul complexes 
 denoted by the same symbols for quadratic algebras and morphisms 
 \cite{YuM2}. Note that $(\mathcal{B}^{!}\otimes \mathcal{C},d)$ is a
cochain $N$-complex of right $\mathcal{C}$-modules, i.e. $d:\mathcal{B}_{n}^{!}\otimes 
\mathcal{C} \rightarrow \mathcal{B}^{!}_{n+1}\otimes \mathcal{C}$ is $\mathcal{C}$-linear. \\
 
 Similarly the Koszul complexes $K(f)$ associated with morphisms 
 $f$ of quadratic algebras generalize as 
 $N$-complexes for morphisms of $N$-homoge\-neous algebras. Let 
 $\calb=A(E,R)$ and $\calc=A(E',R')$ be two $N$-homogeneous algebras 
  and let $f:\calb\rightarrow \calc$ be a morphism of $N$-homogeneous 
  algebras ($f\in \Hom(\calb,\calc)$). One can define the $N$-complex 
 $K(f)=(\calc\otimes \calb^{!\ast}, d)$ by using partial dualization of the $N$-complex $L(f)$ 
 generalizing thereby the construction of \cite{YuM1} or one can 
 define $K(f)$ by generalizing the construction of \cite{YuM2}, 
 \cite{JLL}.\\
  The first way consists in applying the functor
$\mathrm{Hom}_{\mathcal{C}}(-,\mathcal{C})$ to each right $\mathcal{C}$-module of the
$N$-complex $(\mathcal{B}^{!}\otimes \mathcal{C}, d)$. We get a chain $N$-complex of left
$\mathcal{C}$-modules. Since $\mathcal{B}_{n}^{!}$ is a
finite-dimensional vector space,
$\mathrm{Hom}_{\mathcal{C}}(\mathcal{B}_{n}^{!}\otimes
\mathcal{C},\mathcal{C})$ is canonically identified to the left module
$\mathcal{C}\otimes (\mathcal{B}_{n}^{!})^{\ast}$. Then we get the
$N$-complex $K(f)$ whose differential $d$ is easily described in terms
of $f$. In the case $\mathcal{A} = \mathcal{B}= \mathcal{C}$ and
$f=I_{\mathcal{A}}$, this complex will be denoted by $K(\mathcal{A})$.\\ 
 We shall follow hereafter the second more explicit way. 
 Let us associate with $f\in \Hom(\calb, \calc)$ the homogeneous 
 linear mapping of degree zero $\alpha:(\calb^!)^\ast\rightarrow 
 \calc$ defined by setting $\alpha=f:E\rightarrow E'$ in degree 1 
 and $\alpha=0$ in degrees different from 1. The dual $(\calb^!)^\ast$ 
 of $\calb^!$ defined degree by degree is a graded coassociative 
 counital coalgebra and one has $\alpha^{\ast N} =\underbrace{\alpha 
 \ast\dots\ast\alpha}_{N}=0$. Indeed it follows from the definition 
 that $\alpha^{\ast N}$ is trivial in degrees $n\not = N$. On the 
 other hand in degree $N$, $\alpha^{\ast N}$ is the composition
 \[
 R \stackrel{f^{\otimes^N}}{\longrightarrow} E^{\prime\otimes^ 
 N}\longrightarrow E^{\prime\otimes^N}/R'
 \]
 which vanishes since $f^{\otimes^N}(R)\subset R'$. Applying Lemma~1 
 it is easily checked that the $N$-differential
 \[
 d_{\alpha}:\calc \otimes \calb^{!\ast}
 \rightarrow \calc \otimes \calb^{!\ast}
 \]
 coincides with $d$ of the first way.\\
 
 Let us give an even more explicit description of $K(f)$ and pay some 
 attention to the degrees. Recall that by $(\calb^!)^\ast$ we just 
 mean here the direct sum $\oplus_{n}(\calb^!_n)^{\ast}$ of the dual 
 spaces $(\calb^!_n)^{\ast}$ of the finite-dimensional vector 
 spaces $\calb^!_n$. On the other hand, with $\calb=A(E,R)$ as 
 above, one has
 \[
 \calb^!_n=E^{\ast\otimes^n}\ \ \mathrm{if}\ \ n < N
 \]
 and
 \[
 \calb^!_n=E^{\ast\otimes^n}/\sum_{r+s=n-N}E^{\ast\otimes^r}\otimes 
 R^\perp \otimes E^{\ast\otimes^s}\> \> \mathrm{if}\> \>  n\geq N.
 \]
 So one has for the dual spaces
 \begin{equation}
 (\calb^!_n)^{\ast}\cong E^{\otimes^n}\> \>Ê \mathrm{if}\> \>Ê n < N
 \label{eq13}
 \end{equation}
 and
 \begin{equation}
 (\calb^!_n)^{\ast}\cong\bigcap_{r+s=n-N}E^{\otimes^r}\otimes R\otimes 
 E^{\otimes^s}\> \> \mathrm{if}\> \>  n\geq N.
 \label{eq14}
 \end{equation}
 In view of (\ref{eq13}) and (\ref{eq14}), one has canonical injections 
 \[
 (\calb^!_n)^\ast \hookrightarrow (\calb^!_k)^\ast \otimes (\calb^!_\ell)^\ast
\]
for $k+\ell=n$ and one sees that the coproduct $\Delta$ of $(\calb^!)^\ast$ is given by
\[
\Delta(x)=\sum_{k+\ell=n}x_{k\ell}
\]
for $x\in(\calb^!_n)^\ast$ where the $x_{k\ell}$ are the images of $x$ into $(\calb^!_k)^\ast \otimes (\calb^!_\ell)^\ast$ under the above canonical injections.\\
 If $f:\calb\rightarrow \calc = A(E',R')$ is a morphism of 
 $\hg_{N}\alg$, one verifies that the $N$-differential $d$ of $K(f)$ 
 defined above is induced by the linear mappings
 \begin{equation}
c\otimes (e_{1}\otimes e_{2} \otimes \cdots \otimes e_{n}) \mapsto c f(e_{1}) \otimes (e_{2}
\otimes \cdots \otimes e_{n})
 \label{eq15}
 \end{equation}
 of $\calc\otimes E^{\otimes^n} $ into $\calc\otimes E^{\otimes^{n-1}}
 $.
 One has $d\left(\calc_s\otimes (\calb^!_r)^{\ast}\right)\subset 
 \calc_{s+1}\otimes (\calb^!_{r-1})^\ast $ so the $N$-complex $K(f)$ 
 splits into subcomplexes
 \[
 K(f)^n=\oplus_{m}\calc_{n-m}\otimes (\calb^!_m)^{\ast},\ \ \ n\in 
 \mathbb N
 \]
 which are homogeneous for the total degree. Using (\ref{eq13}), 
 (\ref{eq14}), (\ref{eq15}) one can describe $K(f)^0$ as
 \begin{equation}
 \cdots \rightarrow 0 \rightarrow \mathbb K \rightarrow 0 \rightarrow  
\cdots
  \label{eq16}
 \end{equation}
 and $K(f)^n$ as
 \begin{equation}
 \cdots \rightarrow 0\rightarrow 
 E^{\otimes^n}\stackrel{f\otimes I^{\otimes^{n-1}}_{E} }{\longrightarrow} 
E' \otimes E^{\otimes^{n-1}}\rightarrow \cdots \stackrel{I^{\otimes^{n-1}}_{E'}\otimes f 
 }{\longrightarrow} 
 E^{\prime\otimes^n}\rightarrow 0\rightarrow \cdots
 \label{eq17}
 \end{equation}
 for $1\leq n\leq N-1$ while $K(f)^N$ reads
 \begin{equation}
 \cdots 0 \rightarrow R
 \stackrel{f\otimes I^{\otimes^{N-1}}_{E}}{\longrightarrow} 
 E'\otimes E^{\otimes^{N-1}}\rightarrow \cdots \rightarrow E^{\prime\otimes^{N-1}}\otimes E
 \stackrel{can}{\rightarrow} 
 \calc_N\rightarrow 0\cdots
 \label{eq18}
 \end{equation}
where $can$ is the composition of $I^{\otimes^{N-1}}_{E'}\otimes f$ 
with canonical projection of $E^{\prime\otimes^N}$ onto 
$E^{\prime\otimes^N}/R'=\calc_N$. 
\\

Let us seek for conditions of maximal acyclicity for the $N$-complex 
$K(f)$. Firstly, it is clear that $K(f)^0$ is not acyclic, one has 
$_{p}H_0(K(f)^0)=\mathbb K$ for $p\in \{1,\dots,N-1\}$. Secondly 
if $N\geq 3$, it is straightforward that if $n\in\{1,\dots,N-2\}$ 
then $K(f)^n$ is acyclic if and only if $E=E'=0$. Next comes the 
following lemma.

\begin{lemma}
 The $N$-complexes $K(f)^{N-1}$ and $K(f)^N$ are acyclic if and only 
 if $f$ is an isomorphism of $N$-homogeneous algebras.
 \end{lemma}
 
 \noindent\underbar{Proof}. First $K(f)^{N-1}$ is acyclic if and only 
 if $f$ induces an isomorphism $f:E\stackrel{\simeq}{\rightarrow}E'$ of 
 vector spaces as easily verified and then, the acyclicity of $K(f)^N$ 
 is equivalent to $f^{\otimes^N}(R)=R'$ which means that $f$ is an 
 isomorphism of $N$-homogeneous algebras.$\square$\\
 
 It is worth noticing here that for $N\geq 3$ the nonacyclicity of 
 the $K(f)^n$ for $n\in \{ 1,\dots,N-2\}$ whenever $E$ or $E'$ is 
 nontrivial is easy to understand and to possibly cure. Let us 
 assume that $K(f)^{N-1}$ and $K(f)^N$ are acyclic. Then by 
 identifying through the isomorphism $f$ the two $N$-homogeneous 
 algebras, one can assume that $\calb=\calc=\cala=A(E,R)$ and that $f$ 
 is the identity mapping $I_{\cala}$ of $\cala$ onto itself, that is 
 with the previous notation that one is dealing with $K(f)=K(\cala)$. 
 Trying to make $K(\cala)$ as acyclic as possible one is now faced to 
 the following result for $N\geq 3$.
 
 \begin{proposition}
  Assume that $N\geq 3$, then one has 
  \[
  \mathrm{Ker} (d^{N-1}:\cala_2\otimes (\cala^!_{N-1})^\ast\rightarrow 
  \cala_{N+1})=\mathrm{Im} (d:\cala_1\otimes (\cala^!_N)^{\ast}\rightarrow 
 \cala_2\otimes  (\cala^!_{N-1})^\ast)
  \]
  if and only if either $R=E^{\otimes^N}$ or $R=0$.
  \end{proposition}
   
   \noindent \underbar{Proof}. One has 
   \[
   \cala_2\otimes (\cala_{N-1}^!)^\ast
   =E^{\otimes^2}\otimes E^{\otimes^{N-1}}\simeq E^{\otimes^{N+1}},
   \cala_{N+1}\simeq E^{\otimes^{N+1}}/E\otimes R+R\otimes E
   \]
   and 
   $d^{N-1}$ identifies here with the canonical projection 
   \[
   E^{\otimes^{N+1}}\rightarrow E^{\otimes^{N+1}}/E\otimes R+R\otimes E
   \]
   so its kernel is $E\otimes R + R\otimes E$. On the other hand one 
   has $\cala_1\otimes(\cala^!_N)^{\ast}= E\otimes R$ and $d:E\otimes 
   R\rightarrow E^{\otimes^{N+1}}$ is the inclusion. So 
   $\im(d)=\ker(d^{N-1})$ is here equivalent to $R\otimes E=E\otimes 
   R+R\otimes E$ and thus to $R\otimes E=E\otimes R$ since all vector 
   spaces are finite-dimensional. It turns out that this holds if and 
   only if either $R=E^{\otimes^N}$ or $R=0$ (see the appendix).~$\square$
   
   \begin{corollary}
    Assume that $N\geq 3$ and let $\cala=A(E,R)$ be a $N$-homogeneous algebra.
    Then the $K(\cala)^n$ are acyclic for $n\geq N-1$ if and only if 
    either $R=0$ or $R=E^{\otimes^N}$.
    \end{corollary}
    
    \noindent \underbar{Proof}. In view of Proposition 2, $R=0$ or 
    $R=E^{\otimes^N}$ is necessary for the acyclicity of 
    $K(\cala)^{N+1}$; on the other hand if $R=0$ or $R=E^{\otimes^N}$ then 
    the acyclicity of the $K(\cala)^n$ for $n\geq N-1$ is 
    obvious.~$\square$\\
    
    Notice that $R=0$ means that
    $\cala$ is the tensor algebra $T(E)$ 
    whereas $R=E^{\otimes^N}$ means that $\cala=T(E^\ast)^!$. Thus 
    the acyclicity of the $K(\cala)^n$ for $n\geq N-1$ is stable 
    by the duality $\cala\mapsto \cala^!$ as for quadratic algebras 
    ($N=2$). However for $N\geq 3$ this condition does not lead to an 
    interesting class of algebras contrary to what happens for $N=2$ 
    where it characterizes the Koszul algebras \cite{Pri}. This is the very 
    reason why another generalization of Koszulity has been 
    introduced and studied in \cite{RB3} for $N$-homogeneous algebras.
  
 \section{Koszul homogeneous algebras}

Let us examine more closely the $N$-complex $K(\mathcal{A})$:
$$\cdots \longrightarrow \mathcal{A}\otimes (\cala^!_i)^\ast \stackrel{d}{\longrightarrow} \mathcal{A}\otimes (\cala^!_{i-1})^\ast \longrightarrow \cdots
\longrightarrow \mathcal{A}\otimes (\cala^!_1)^\ast
\stackrel{d}{\longrightarrow} \mathcal{A} \longrightarrow 0\,.$$
The $\cala$-linear map $d:\cala\otimes (\cala^!_i)^\ast\rightarrow \cala\otimes (\cala^!_{i-1})^\ast$ is induced by the canonical injection (see in last section)
\[
(\cala^!_i)^\ast\hookrightarrow (\cala^!_1)^\ast \otimes (\cala^!_{i-1})^\ast=\cala_1\otimes (\cala^!_{i-1})^\ast \subset \cala\otimes (\cala^!_{i-1})^\ast.
\]
The degree $i$ of
$K(\mathcal{A})$ as $N$-complex has not to be confused with the total
degree $n$. Recall that, when $N=2$, the quadratic algebra $\mathcal{A}$ is said to be Koszul if
$K(\mathcal{A})$ is acyclic at any degree $i>0$ (clearly it is equivalent to
saying that each complex $K(\mathcal{A})^{n}$ is acyclic for any total degree
$n>0$). 

\noindent For any $N$, it is possible to contract the
$N$-complex $K(\mathcal{A})$ into (2-)complexes by putting
together alternately $p$ or $N-p$ arrows $d$ in $K(\mathcal{A})$. The
complexes so obtained are the following ones
$$\cdots
\stackrel{d^{N-p}}{\longrightarrow} \mathcal{A}\otimes (\cala^!_{N+r})^\ast
\stackrel{d^{p}}{\longrightarrow} \mathcal{A}\otimes (\cala^!_{N-p+r})^\ast
\stackrel{d^{N-p}}{\longrightarrow} \mathcal{A}\otimes (\cala^!_{r})^\ast
\stackrel{d^{p}}{\longrightarrow} 0\,,$$
which are denoted by $C_{p,r}$. All the possibilities are covered by
the conditions  $0 \leq r \leq N-2$ and $r+1 \leq p \leq N-1$. Note
that the complex $C_{p,r}$ at degree $i$ is $\mathcal{A}\otimes
(\cala^!_k)^\ast$, where $k=jN+r$ or $k=(j+1)N-p+r$, according to
$i=2j$ or $i=2j+1$ ($j \in \mathbb N$).

In \cite{RB3}, the complex $C_{N-1,0}$ is called the \emph{Koszul complex} of
$\mathcal{A}$, and the homogeneous algebra $\mathcal{A}$ is said to be \emph{Koszul} if this complex
is acyclic at any degree $i>0$. A motivation for this definition is
that Koszul property is equivalent to a purity property of the
minimal projective resolution of the trivial module. One has the following result \cite{RB3}, \cite{RB4} : 
\begin{proposition}
Let $\mathcal{A}$ be a homogeneous algebra of degree $N$. For $i=2j$ or $i=2j+1$, $j \in \mathbb N$, the graded vector space
$\mathrm{Tor}_{i}^{\mathcal{A}}(\mathbb K, \mathbb K)$ lives in
degrees $\geq jN$ or $\geq jN+1$ respectively. Moreover, $\mathcal{A}$ is Koszul if and
only if each $\mathrm{Tor}_{i}^{\mathcal{A}}(\mathbb K, \mathbb K)$
is concentrated in degree $jN$ or $jN+1$ respectively (purity property).
\end{proposition}
When $N=2$, it is exactly Priddy's definition \cite{Pri}. Another
motivation is that a certain cubic Artin-Schelter regular algebra has the
purity property, and this cubic algebra is a good candidate for making
non-commutative algebraic geometry \cite{AS}, \cite{ATVB}. Some other non-trivial
examples are contained in \cite{RB3}.

The following result shows how the Koszul complex $C_{N-1,0}$ plays a
particular role. Actually all the other contracted complexes of
$K(\mathcal{A})$ are irrelevant as far as acyclicity is concerned. 
\begin{proposition}
Let $\mathcal{A}=A(E,R)$ be a homogeneous algebra of degree $N\geq 3$. Assume
that $(p,r)$ is distinct from $(N-1,0)$ and that $C_{p,r}$ is exact at
degree $i=1$. Then $R=0$ or $R=E^{\otimes ^{N}}$.
\end{proposition}
\noindent \underbar{Proof}. 
Assume $r=0$, hence $1\leq p \leq N-2$. Regarding $C_{p,0}$ at degree
1 and total
degree $N+1$, one gets the exact sequence
$$E\otimes R
\stackrel{d^{p}}{\longrightarrow} E^{\otimes ^{N+1}}
\stackrel{d^{N-p}}{\longrightarrow} E^{\otimes ^{N+1}}/ E\otimes R + R
\otimes E,$$
where the maps are the canonical ones. Thus $E\otimes R =E\otimes R + R
\otimes E$, leading to $R\otimes E= E\otimes R$. This holds only if
$R=0$ or $R=E^{\otimes ^{N}}$ (Appendix).\\

Assume now $1\leq r \leq N-2$ (hence $r+1\leq p \leq N-1$). Regarding
$C_{p,r}$ at degree 1 and total
degree $N+r$, one gets the exact sequence
$$(\cala^!_{N+r})^\ast
\stackrel{d^{p}}{\longrightarrow} E^{\otimes ^{N+r}}
\stackrel{d^{N-p}}{\longrightarrow} E^{\otimes ^{N+r}}/ R\otimes
E^{\otimes ^{r}},$$
where the maps are the canonical ones. Thus $(\cala^!_{N+r})^\ast=R\otimes
E^{\otimes ^{r}}$, and $R\otimes E^{\otimes ^{r}}$ is contained in
$E^{\otimes ^{r}}\otimes R$. So $R\otimes E^{\otimes ^{r}}=E^{\otimes
  ^{r}}\otimes R$, which implies again $R=0$ or $R=E^{\otimes ^{N}}$ (Appendix).$\square$\\

It is easy to check that, if $R=0$ or $R=E^{\otimes ^{N}}$, any
$C_{p,r}$ is exact at any degree $i>0$. On the other hand, for any $R$, one has
$$H_{0}(C_{p,r}) = \oplus_{0\leq j \leq N-p-1}\> \>  E^{\otimes ^{j}}
\otimes E^{\otimes ^{r}},$$
which can be considered as a Koszul left $\mathcal{A}$-module if $\cala$ is Koszul. 

\section{Appendix : a lemma on tensor products}
\begin{lemma}
Let $E$ be a finite-dimensional vector space. Let $R$ be a subspace of
$E^{\otimes ^{N}}$, $N\geq 1$. If $R\otimes E^{\otimes ^{r}}=E^{\otimes
  ^{r}}\otimes R$ holds for an integer $r\geq 1$, then $R=0$ or
$R=E^{\otimes ^{N}}$.
\end{lemma}
\noindent\underbar{Proof}. 
Fix a basis $X=(x_{1}, \ldots ,x_{n})$ of $E$, ordered by
$x_{1}< \cdots <x_{n}$. The set $X^{N}$ of the words of length $N$ in
the letters $x_{1}, \ldots ,x_{n}$ is a basis of $E^{\otimes ^{N}}$
which is lexicographically ordered. Denote by $S$ the
$X^{N}$-reduction operator of $E^{\otimes ^{N}}$ associated to
$R$ \cite{RB1}, \cite{RB2}. This means the following properties:
\\ \\
(i) $S$ is an endomorphism of the vector space $E^{\otimes ^{N}}$ such
that $S^{2}=S$,
\\ \\
(ii) for any $a\in X^{N}$, either $S(a)=a$ or $S(a)<a$ (the
latter inequality means $S(a)=0$, or otherwise any word occuring in
the linear decomposition of $S(a)$ on $X^{N}$ is $<a$ for the
lexicographic ordering),
\\ \\
(iii) $\mathrm{Ker}(S)=R$.
\\ \\
Then $S\otimes I_{E^{\otimes ^{r}}}$ and $I_{E^{\otimes ^{r}}} \otimes
S$ are the $X^{N+r}$-reduction operators of $E^{\otimes ^{N+r}}$,
respectively associated to $R\otimes E^{\otimes ^{r}}$ and $E^{\otimes^{r}}\otimes R$. By assumption these endomorphisms are equal. In particular, one has
$$\mathrm{Im}(S)\otimes E^{\otimes ^{r}}= E^{\otimes ^{r}}\otimes
\mathrm{Im}(S).$$
But the subspace $\mathrm{Im}(S)$ is monomial, i.e. generated
by words. So it suffices to prove the lemma when $R$ is monomial.

Assume that $R$ contains the word
$x_{i_{1}}\ldots x_{i_{N}}$. For any letters $x_{j_{1}}, \ldots ,
x_{j_{r}}$, the word $x_{i_{1}}\ldots x_{i_{N}}x_{j_{1}} \ldots
x_{j_{r}}$ belongs to $E^{\otimes ^{r}}\otimes R$. Thus
$x_{i_{r+1}}\ldots x_{i_{N}}x_{j_{1}} \ldots x_{j_{r}}$ belongs to
$R$. Continuing the process, we see that any word belongs to $R$.~$\square$

\end{document}